\documentclass[11pt]{article}
\usepackage{amsmath,amssymb}
\usepackage{eucal}
\title{\leavevmode\vadjust{\vskip -5mm}
\textbf{Canonical Deformed Groups of Diffeomorphisms and\\
Finite Parallel Transports in Riemannian Spaces}}
\author{\sc Serhiy E. SAMOKHVALOV\thanks{{\bf e}-{\it mail}: samokhval@dstu.dp.ua}\\
 \tabaddress{Department of Applied Mathematics \\ State Technical University, Dniprodzerzhinsk,
 Ukraine}}

\date{April 22, 2007}
\def\tabaddress#1{{\small\it\begin{tabular}[t]{c}#1
\\[1.2ex]\end{tabular}}}

\def\<#1>{\langle#1\rangle}

\def\beq{\begin{equation}}
\def\eeq{\end{equation}}
\def\bea{\begin{eqnarray}}
\def\eea{\end{eqnarray}}
\def\beann{\begin{eqnarray*}}
\def\eeann{\end{eqnarray*}}
\def\ben{\begin{enumerate}}
\def\een{\end{enumerate}}

\def\qed{\ifvmode\removelastskip\fi
{\unskip\nobreak\hfil\penalty50\hbox{}\nobreak\hfil \hbox{\vrule
height1.2ex width1.2ex}\parfillskip=0pt \finalhyphendemerits=0
\par\smallskip}}

\def\texthook{\vrule height 0pt depth 0.4pt width 3.5pt
          \vrule height 5pt depth 0.4pt \kern 3pt}
\def\scripthook{\vrule height 0pt depth 0.2pt width 1.5pt
                \vrule height 3pt depth 0.2pt width 0.2pt \kern 1pt}

\begin{document}

\maketitle

\thispagestyle{empty}

\begin{abstract}

\noindent We show that finite parallel transports of vectors in
Riemannian spaces, determined by the multiplication law in the
deformed groups of diffeomorphisms, and sequences of infinitesimal
parallel transports of vectors along geodesics are equivalent.

\end{abstract}


\\

\clearpage


Performing group-theoretic description of spaces with torsion-free
affine connection and Riemannian spaces with the help of deformed
groups of diffeomorphisms $\Gamma^H_T$, in work [1] we introduced
the concept of parallel transport of vectors to finite and not
just infinitesimal distance $\tilde t=x'-x$, the concept following
from multiplication law in groups $\Gamma^H_T$ . On the other
hand, in classic approach finite parallel transports are generated
from infinitesimal ones in case of movement along curves [2], in
curved space the result depending upon the selected curve. In case
of infinitesimal shifts $\tilde t$  both transports yield
identical results, affine vielbein field and affine connection
coefficients, defined by the action of deformed group of
diffeomorphisms $\Gamma^H_T$ , are defined by the first two orders
in expansion of functions of deformations with respect to $\tilde
t$, while all higher orders in spaces with affine connection
remain undefined though influence the result of finite parallel
transports.

In this work we show that requirement of that for two arbitrary
points $x$ and $x'$ in arbitrarily curved space there exists at
least one curve connecting them, for which sequence of
infinitesimal parallel transports would yield result identical to
finite parallel transport (specified by multiplication law of
deformed group of diffeomorphisms $\Gamma^H_T$ ) is enough to
prove that such curve is (locally) unique and is a geodesic. In
this approach the derivatives from the deformation functions in
the direction of geodesic are the first integrals of system of
differential equations which specify it. By the latter deformation
functions (and not just their first two orders in expansion with
respect to shifts $\tilde t$ ) are uniquely defined with respect
to affine vielbein field and affine connection coefficients, thus
finite parallel transport to the distance $\tilde t=x'-x$
specified by multiplication law of group $\Gamma^H_T$  is unique
and yields the result of parallel transport along the geodesic
connecting points $x$  and $x'$ .

Groups $\Gamma^H_T$ produced by means of such deformations are a
certain generalization of finite parametric canonical Lie groups
[3] for infinite case and, therefore, are also called canonical as
well as to deformations with the help of which they are built.

In the Riemannian space deformation functions (including their
second order in expansion with respect to the shifts specified by
connection) are completely defined with respect to the orthonormal
vielbein field by requirement for vectors not to change their
length and only rotate in case of finite parallel transports in
group-theoretic mode. This work shows that such deformations are
canonical.

We also show that in Riemannian space, where geodesics in natural
parameterization are extremals of energy functional which is
defined by metric, the function of central extremal fields, which
start from each point [4] (action function) uniquely (with
accuracy to the choice of orthonormal vielbein field coordinated
with metric) define deformation functions. Moreover,
Hamilton-Jacobi equation, which is satisfied by action function,
for deformation function comes down to the equation which follows
from the requirement of invariability of vector length in case of
finite parallel transports.

The determined relation between the first integrals of geodesics
as well as action function with deformation functions results in
their new understanding and can be of applied significance,
particularly in gravitation theory.

This article uses the same notations and assumptions as [1].
Specifically, all relations are obtained with the bounds of single
coordinate chart with fixed though arbitrary coordinates.

\section{Deformed groups of diffeomorphisms}

Let in coordinate chart $O$ with coordinates  $x^\mu$ (coordinate
indices are selected from Greek alphabet) there act a local group
of diffeomorphisms in additive parameterization $\Gamma_T=\{\tilde
t(x)\}$ (\emph{undeformed group of diffeomorphisms}) [1] according
to the formula
\[x'^\mu = x^\mu+ \tilde t^\mu(x).
\]
Smooth functions $ \tilde t^\mu(x)$  parameterizing group
$\Gamma_T$ satisfy the requirement $\{\delta^\mu_\nu+\partial_\nu
\tilde t^\mu(x)\}\neq0$, $ \forall x \in O $, where
$\partial_\nu:=\partial/\partial x^\nu$, and multiplication law in
it $\tilde t''= \tilde t\times \tilde t'$ is specified by the
formula:
\begin{equation}\label{eq1}
\tilde t''^\mu (x)= \tilde t^\mu (x)+ \tilde t'^\mu (x').
\end{equation}

\emph{Deformed group of diffeomorphisms}  $\Gamma^H_T = \{t(x)\}$
[1] is parameterized by the functions $t^m(x)$  (we will use
indices from Latin alphabet for them) is produced from group
$\Gamma_T$ as isomorphism specified by deformation $H$ according
to the formula
\begin{equation}\label{eq2}
t^m (x)=H^m (x, \tilde t (x))
\end{equation}
\noindent with the help of smooth deformation functions $H^m(x,
\tilde t(x))$ with properties

$1H)\ H^m(x, 0)=0,\ \forall x \in O$;

$2H)\ \exists$ functions  $K^\mu(x,t(x)):K^\mu (x,H(x, \tilde
t(x)))= \tilde t^\mu(x),\ \forall x \in O,\ \tilde t \in
\Gamma_T$.

Functions $t^m(x)$ satisfy the condition det $\{\delta^\mu_\nu +
d_\nu K^\mu (x,t(x))\}\neq 0,\ \forall x\in O$, where
$d_\nu:=d/dx^\nu$. Functions $K^\mu(x,t(x))$ which are present in
property $2H)$ specify reverse transition $\Gamma^H_T \rightarrow
\Gamma_T$ . Multiplication law $t''=t \ast t'$ in group
$\Gamma^H_T$ is defined by multiplication law (1) in group
$\Gamma_T$ and isomorphism (2):
\begin{equation}\label{eq3}
t''^m (x)=\varphi^m (x,t(x), t'(x')):=H^m
(x,K(x,t(x)+K(x',t'(x'))),
\end{equation}
\noindent where
\begin{equation}\label{eq4}
x'^\mu = f^\mu (x,t(x)):=x^\mu + K^\mu (x,t(x)).
\end{equation}
\noindent Group  $\Gamma^H_T$ acts smoothly in chart $O$ according
to the formula (4).

With the help of functions $\varphi(x,t,t')$  which specify
multiplication law (3) auxiliary matrices are defined:
\begin{equation}\label{eq5}
\lambda (x,t)^m_n:= \partial_{n'} \varphi^m (x,t,t')|_{t'=0} =
h(x+ \tilde t)^\mu_n \partial_{\tilde \mu} H^m (x, \tilde t)|_
{\tilde t=K(x,t)},
\end{equation}
\begin{equation}\label{eq6}
\mu(x,t)^m_n:= \partial_{n'} \varphi^m (x,t',t)|_{t'=0} =
h(x)^\nu_n (\delta^\mu_\nu+\partial_\nu \tilde
t^\mu)\partial_{\tilde \mu} H^m(x, \tilde t)|_{\tilde t=K(x,t)},
\end{equation}
\noindent where  $h(x)^\mu_m :=\partial_m K^\mu(x,t)|_{t=0}$ and
$\partial_m:= \partial / \partial t^m$ (primed index standing for
differentiation with respect to $t'$ and tilded one for that with
respect to $\tilde t$ ).

Let's emphasize that despite the fact that deformations are
defined as isomorphisms of diffeomorphisms group, in general case
they cannot be compensated by coordinate transformations in chart
$O$ which result only in internal automorphisms of diffeomorphisms
group.

As shown in work [1], deformed group of diffeomorphisms
$\Gamma^H_T$ specifies by its action in chart $O$ geometric
structure of space with torsion-free affine connection and
deformations lead to change of its characteristics, specifically
to curvature.

Generators  $X_m=h(x)^\mu_m \partial_\mu$ of action (4) of group
$\Gamma^H_T$ specify in $O$ an affine vielbein field, matrices
$h(x)^\mu_m$ and matrices $h(x)^m_\mu$ reverse to them transit
between coordinate and affine veilbeins.

Element of group  $\Gamma^H_T$ are accorded with vector fields
$t=t^m(x)X_m$ and parameters $t^m(x)$ of group $\Gamma^H_T$ appear
for components of these fields in basis $X_m$ .

Multiplication law (3) in group  $\Gamma^H_T$ written down for
elements  $t$ and $\tau$ for the case of infinitesimal second
multiplier
\[ (t \ast \tau)^m(x)=t^m(x)+\lambda (x,t(x))^m_n \tau^n(x'),
\]
\noindent defines the rule of composition of vectors fitted in
different points or the \emph{rule of parallel transport} of
vector field $\tau$ from point $x'$ to point $x$ for finite
distance $x'-x= \tilde t=K(x,t)$:
\[ \tau^m_\parallel (x)=\lambda(x,t)^m_n\tau^n(x').
\]
\noindent In coordinate basis this relation appears as
\begin{equation}\label{eq7}
\tau^\mu_\parallel (x)=\lambda(x,\tilde
t)^\mu_\nu\tau^\nu(x')=\partial_{\tilde \nu}H^\mu(x, \tilde
t)\tau^\nu(x'),
\end{equation}
\noindent (here, in particular, $H^\mu(x,\tilde
t)=h(x)^\mu_mH^m(x,\tilde t)$ ) or in case of infinitesimal
$\tilde t$:
\[ \tau^\mu_\parallel (x)=\tau^\mu(x)+\tilde
t^\nu(x)\nabla_\nu\tau^\mu(x),
\]
\noindent where
\[
\nabla_\nu\tau^\mu(x)=\partial_\nu\tau^\mu(x)+\Gamma(x)^\mu_{\sigma\nu}\tau^\sigma(x)
\]
\noindent is a covariant derivative while functions
\begin{equation}\label{eq8}
\Gamma(x)^\mu_{\sigma\nu}:=\partial^2_{\tilde \sigma \tilde
\nu}H^\mu(x,\tilde t)|_{\tilde t=0}
\end{equation}
\noindent acquire geometric meaning of (torsion-free) affine
connection coefficients in coordinate basis. The higher orders of
expansion of function  $H^\mu(x,\tilde t)$ with respect to $\tilde
t $ do not influence the connection though, as follows from
formula (7), define the result of finite parallel transports thus
require precise definition which is given in the next section of
the work. The skew-symmetric part
\begin{equation}\label{eq9}
R^m_{lkn}:=\rho^m_{lkn}-\rho^m_{lnk}
\end{equation}
\noindent of coefficients
$\rho^m_{lkn}:=\partial^2_{lk}\mu(x,t)^m_n|_{t=0}$ (here and later
dependence upon $x$ is not explicitly shown for cases where it is
obvious) written down in coordinate basis appears as [1]:
\begin{equation}\label{eq10}
R^\mu_{\lambda k
\nu}=\partial_k\Gamma^\mu_{\nu\lambda}-\partial_\nu\Gamma^\mu_{k\lambda}+
\Gamma^\mu_{k\sigma}\Gamma^\sigma_{\nu\lambda}-\Gamma^\mu_{\nu\sigma}\Gamma^\sigma_{k\lambda},
\end{equation}
\noindent consequently it is Riemann-Cristoffel curvature tensor.
In Riemannian space deformation functions are satisfied by
equations
\begin{equation}\label{eq11}
\partial_{\tilde \mu}H^\rho(x,\tilde t)\partial_{\tilde \nu}H^\sigma(x, \tilde
t)g(x)_{\rho\sigma}=g(x+\tilde t)_{\mu\nu},
\end{equation}
\noindent which follows from the requirement for vectors not to
change their length and just rotate in case of finite parallel
transports. Equation (11) allows to uniquely define deformation
functions in coordinate basis  $H^\mu(x,\tilde t)$ according to
metric, specifically their second order in expansion with respect
to $\tilde t$, i.e. according to (8) affine connection
coefficients which become equal to Cristoffel symbols:
\begin{equation}\label{eq12}
\Gamma^\rho_{\mu\nu}=\frac{1}{2}g^{\rho\sigma}(\partial_\mu
g_{\nu\sigma}+\partial_\nu g_{\mu\sigma}-\partial_\sigma
g_{\mu\nu}).
\end{equation}

\section{Canonical deformations}

Let's assume that group $\Gamma^H_T$  is such that between two
arbitrary points $x$ , $x'\in O$  there exists parametric curve
$x(s),\ s \in I=[0,T]:\ x(0)=x,\ x(T)=x'$ , along which all
tangent vectors $\tau(s):=\dot x(s)$ (point stands for
differentiation with respect to $s$) parallel to each other in
terms of parallel finite transports, this is according to (7)
\begin{equation}\label{eq13}
\tau^\mu(s)=\lambda(x(s), \tilde t(s,s'))^\mu_\nu\tau^\nu(s'),\
\forall s,\ s'\in I,
\end{equation}
\noindent where $\tilde t(s,s'):=x(s')-x(s)$. Condition (13) for
deformation function becomes:
\begin{equation}\label{eq14}
\tau^\mu(s)=\partial_{\tilde \nu}H^\mu(x(s), \tilde
t(s,s'))\tau^\nu(s')=\frac{d}{ds'} H^\mu(x(s), \tilde t(s,s')),\
\forall s,\ s' \in I.
\end{equation}

\textbf{\ Definition.} \emph{Deformations are called
\textbf{canonical} and groups $\Gamma^H_T$  produces with their
help are called \textbf{canonical deformed groups of
diffeomorphisms}, provided that between two arbitrary points $x,\
x' \in O$, there exists at least one curve  $x(s)$, along which
deformation functions satisfy the requirement (14)}.

Differentiation of equation (14) with respect to $s'$  with
$s'=s$, with consideration of formula (8) yields the equation of
curve  $x(s)$, the existence of which is assumed in Definition:
\begin{equation}\label{eq15}
\dot
\tau^\mu(s)+\Gamma(x(s))^\mu_{\sigma\nu}\tau^\sigma(s)\tau^\nu(s)=0,\
\dot x^\mu(s)=\tau^\mu(s),
\end{equation}
\noindent which is evidence of the fact that curve $x(s)$  is
geodesic in affine parameterzation for the structure of affine
connection specified in $O$ by the action of group  $\Gamma^H_T$.
Locally between two points $x$  and $x'$  only one geodesic can be
drawn, accordingly \emph{the requirement of existence of at least
one curve for which condition (14) is satisfied is enough for
identification that such curve is geodesic}.

Now let's rewrite conditions (13) and (14) for $s=0$  (and
substituting $s'$ for $s$ ) with consideration of definition of
shift $\tilde t(0,s)=x(s)-x$:
\begin{equation}\label{eq16}
\tau^\mu=\lambda(x,x'-x)^\mu_\nu\tau^{'\nu}=H^\mu(x,x'-x)=:u^\mu(x,x',\tau').
\end{equation}
\noindent Here we assumed $\tau:=\tau(0)$  as well as denoted
current points and tangent vectors to curve as  $x':=x(s), \
\tau':=\tau(s)$. Thus functions $u^\mu(x,x',\tau')$ , being
derivatives from deformation functions along geodesics, are
constant for canonical deformation and are independent first
integrals of autonomous system of differential equations of
geodesics (15). Thus, functions  $u^\mu(x,x',\tau')$ are to
satisfy the equation [5]
\begin{equation}\label{eq17}
\frac{\partial u^\mu}{\partial x^{'\nu}}\tau^{'\nu}-\frac{\partial
u^\mu}{\partial
\tau^{'\nu}}\Gamma(x')^\nu_{\rho\sigma}\tau^{'\rho}
\tau^{'\sigma}=0,
\end{equation}
\noindent geodesics being characteristics of which. With
consideration of relation (16) equation (17) comes down to the
following equation for auxiliary matrices $\lambda(x, \tilde
t)^\mu_\nu$ :
\begin{equation}\label{eq18}
[\partial_{\tilde \rho}\lambda(x, \tilde t)^\mu_\sigma - \Gamma(x+
\tilde t)^\nu_{\rho\sigma}\lambda(x, \tilde
t)^\mu_\nu]\tau^{'\rho}\tau^{'\sigma}=0,
\end{equation}
\noindent which in terms of deformation functions may become:
\begin{equation}\label{eq19}
[\partial^2_{\tilde \rho \tilde \sigma}H^\mu(x, \tilde t) -
\Gamma(x+ \tilde t)^\nu_{\rho\sigma}\partial_{\tilde \nu}H^\mu(x,
\tilde t)]\tau^{'\rho}\tau^{'\sigma}=0.
\end{equation}

Boundary conditions for auxiliary matrices  $\lambda(x, \tilde
t)^\mu_\nu$ and deformation functions $H^\mu(x, \tilde t)$  are
specified in point $x$ and follow from properties 1H), 2H):
\begin{equation}\label{eq20}
\lambda(x, 0)^\mu_\nu=\delta^\mu_\nu,
\end{equation}
\begin{equation}\label{eq21}
H^\mu(x, 0)=0, \ \partial_\nu H^\mu(x, 0)=\delta^\mu_\nu,
\end{equation}
\noindent which in selection of initial vector $\tau$ defines
boundary conditions for the first integrals:
\begin{equation}\label{eq22}
u^\mu(x,x,\tau)=\tau^\mu.
\end{equation}
\noindent According to the theorem of existence and uniqueness of
solution to linear differential equations in partial derivatives
[5], Cauchy problem (17), (22) for arbitrary symmetric with
respect to the lower indices and smooth functions
$\Gamma(x)^\nu_{\rho\sigma}$ has unique solution which corresponds
to geodesic starts from the point $x$ in the direction of vector
$\tau$ . This provides for unique solvability of problem (18),
(20) for auxiliary matrices $\lambda(x, \tilde t)^\mu_\nu$ , as
well as problem (19), (20) for deformation functions in coordinate
basis $H^\mu(x, \tilde t)$ which, accordingly, correspond to the
central field of geodesics which start from point $x$ with
different initial vectors $\tau$. \emph{Thus, arbitrary
(torsion-free) affine connection coefficients
$\Gamma(x)^\nu_{\rho\sigma}$ for canonical deformations uniquely
define deformation functions in coordinate basis}.

The zeroth and first orders of expansion of function  $H^\mu(x,
\tilde t)$ with respect to shifts $\tilde t$  are defined by
boundary conditions (21). Sequentially differentiating equation
(19) with respect to $s$ in zero we derive the following recurrent
formula for finding coefficients of expansion of functions
$H^\mu(x, \tilde t)$  with respect to $\tilde t$ in any order
$n\geq2$ :
\[ \partial^{(n)}_{\tilde \nu_1,..., \tilde \nu_n}H^\mu(x,0)=
\partial^{(n-2)}_{\{\tilde \nu_1,..., \tilde
\nu_{n-2}}[\Gamma(x+\tilde t)^\sigma_{\tilde \nu_{n-1} \tilde
\nu_n\}}\partial_{\tilde \sigma}H^\mu(x, \tilde t)]|_{\tilde t=0},
\]
\noindent where curly brackets, as usual, stand for symmetrization
with respect to indices placed in them. Specifically, with $n=2$
we derive relation (8) and with $n=3$ :
\[
\triangle^\mu_{\nu\rho\tau}=\partial_{\{\nu}\Gamma^\mu_{\rho\tau\}}+
\Gamma^\sigma_{\{\nu\rho}\Gamma^\mu_{\tau\}\sigma},
\]
\noindent where $\triangle^\mu_{\nu\rho\tau}:=\partial^3_{\tilde
\nu \tilde \rho \tilde \tau}H^\mu(x, \tilde t)|_{\tilde t=0}$ ,
wherefrom, as shown in [1], it follows that coefficients
$\rho^\mu_{\nu\rho\tau}$ are defined by curvature tensor (10)
uniquely:
\begin{equation}\label{eq23}
\rho^\mu_{\nu\rho\tau}=
\frac{1}{3}(R^\mu_{\nu\rho\tau}+R^\mu_{\rho\nu\tau}).
\end{equation}
\noindent and definition (9) comes down to the known identity
$R^\mu_{\nu\rho\tau}+R^\mu_{\rho\tau\nu}+R^\mu_{\tau\nu\rho}=0$.

Let's summarize the derivations.

\textbf{Proposition 1 (of uniqueness of canonical deformations).}

\emph{a) The condition of canonicity of deformations (14) locally
uniquely defines curve between points $x$ and $x'$ , along which
it is performed, this curve being geodesic in affine
parameterization of space with (torsion-free) affine connection,
the structure of which is specified in $O$ by the action of group
$\Gamma^H_T$ .}

\emph{b) The functions of canonical deformations in coordinate
basis $H^\mu(x, \tilde t)$ are uniquely defined by coefficients
$\Gamma(x)^\sigma_{\rho\sigma}$ of (torsion-free) affine
connection in a single-valued manner, and derivatives from them
along geodesics are the first integrals of system of differential
equations for geodesics.}

It is easy to verify that for three arbitrary points $x,\ x'$ ,
and $x''$ which lie on one geodesic, matrix
$\lambda(x,x'-x)^\mu_\sigma\lambda(x',x''-x')^\sigma_\nu$ do not
depend upon $x'$ and also satisfy the equation (18) and the
condition (20). From uniqueness of solution to problem (18), (20),
this is provided for by composition law of parallel transport of
arbitrary vectors $\theta$ along geodesic:
\begin{equation}\label{eq24}
\theta^\mu(x)=\lambda(x,x''-x)^\mu_\nu\theta^\nu(x'')=
\lambda(x,x'-x)^\mu_\sigma\lambda(x',x''-x')^\sigma_\nu\theta^\nu(x''),
\end{equation}
\noindent which, in its turn, leads to coincidence of
\emph{arbitrary sequence of parallel vector transports along
geodesic (including integral sequence of infinitesimal transports)
with resulting finite transport.} And vice versa, direct
verification shows that if composition law (24) is performed for
three arbitrary points  $x,\ x'$, and $x''$ on geodesic and vector
$\theta$  for finite parallel transports, the functions
$\lambda(x,x'-x)^\mu_\nu$, with the help of which the law is
performed, satisfy the equation (18), accordingly, the
\emph{respective deformation will be canonic.}

\textbf{Proposition 2 (of finite parallel transports).}

\emph{For finite parallel transport of arbitrary vectors $\theta$
to the distance  $\tilde t=x'-x$ to yield the result of integral
sequence of infinitesimal parallel transports along geodesic
connecting points $x$ and $x'$ , it is necessary and sufficient
for deformed group of diffeomorphisms $\Gamma^H_T$ , defining the
transport, to be canonic. }

Let's point out that the condition of canonicity of deformations
(14) is applied to the functions of deformation in coordinate
vielbein  $H^\mu(x, \tilde t)$, accordingly the selection of
affine veilbein $X_m$ remains arbitrary and does not influence
upon the canonicity of deformation.

That is, if deformation with deformation functions  $H^m(x, \tilde
t)$ is canonic, canonic is the deformation with functions
\[ H^{'m}(x, \tilde t)=L(x)^m_nH^n(x, \tilde t)
\]
\noindent where $L(x)^m_n$ are arbitrary nondegenerated matrices
dependent upon $x$.

\section{Criteria of canonicity}

Let's study the parameters of deformed group of diffeomorphisms
$\Gamma^H_T$ which correspond to shifts $\tilde
t(s,s'):=x(s')-x(s)$ along curve $x(s)$:
\begin{equation}\label{eq25}
t^m(s,s')=H^m(x(s), \tilde t(s,s')).
\end{equation}
\noindent In their terms the condition of deformations canonicity
(14) may become
\[ \tau^m(s):=\frac{d}{ds'}t^m(s,s'),\ \forall s, \ s' \in I,
\]
\noindent where $\tau^m(s)=h(x(s))^m_\mu\tau^\mu(s)$  are
components of vector in affine basis $X_m$ which is tangent to
curve in point $x(s)$, wherefrom follows condition for parameters
of group  $\Gamma^H_T$:
\begin{equation}\label{eq26}
t^m(s,s')=(s'-s)\tau^m(s),\ \forall s, \ s' \in I,
\end{equation}
\noindent equivalent to condition (14). Turning equation (25) with
the use of property 2 H), this condition can have the following
equivalent appearance:
\begin{equation}\label{eq27}
\tilde t^\mu(s,s')=x(s')^\mu-x(s)^\mu=K^\mu(x(s),(s'-s)\tau(s))
\quad \forall s, \ s' \in I,
\end{equation}
\noindent Thus, the curve for which condition (27) is fulfilled,
according to the first part of Proposition 1 is geodesic passing
through point  $x(s)$ in direction of vector  $\tau(s)$.

The opposite is true.

\textbf{Proposition 3 (criterion of canonicity for geodesic
equation).}

\emph{If arbitrary geodesic in affine parameterization in chart
$O$ is written down as (27) with functions $K^\mu(x,t)$ , defined
by the deformation functions with property 2H), such deformations
will be canonic.}

This can be easily verified directly by fulfilling double
differentiation of equation (27) with respect to  $s'$ (assuming
$s=0$ and redesignating  $s'$ as $s$ ) and demanding for
coincidence of derivations with geodesic equation (15), which
leads to equation (19) for deformation functions, providing for
their canonicity.

Let's examine the product of elements of canonic deformed group
$\Gamma^H_T$ which correspond to shifts along geodesic and,
therefore, have representation (26):
\begin{equation}\label{eq28}
(s_1+s_2)\tau^m(s)=\varphi^m(x(s),s_1\tau(s),s_2\tau(s')),
\end{equation}
\noindent where  $s'=s+s_1$. This formula generalizes
\emph{canonic multiplication law} [3] for canonic finite parameter
Lie groups for the case of infinite groups $\Gamma^H_T$ , which
accounts for our customary term "canonic" for groups $\Gamma^H_T$
and respective deformations. Formula (28) defines homomorphism of
the additive Abelian group $T^1=\{s\}$ into group $\Gamma^H_T$ .

Differentiating (28) with respect to $s_{2}$ in zero we derive
\begin{equation}\label{eq29}
\tau^m(s)=\lambda(x(s),s_1\tau(s))^m_n\tau^n(s'),
\end{equation}
\noindent which with consideration of formula (5) is equivalent to
condition of canonicity of deformations (14).

Let's study the "left" analogue of condition (29). To this end
let's differentiate (28) with respect to  $s_1$ in zero. We derive
as the result:
\begin{equation}\label{eq30}
\tau^m(s)=\mu(x(s),s_2\tau(s))^m_n\tau^n(s)+s_2\tau^m(s),
\end{equation}
Both equations (29) and (30) fix the curve along which canonic
multiplication law (28) is fulfilled. Thus, differentiating in
zero (29) with respect to $s_1$ , or (30) with respect to  $s_2$,
we derive the equation:
\begin{equation}\label{eq31}
\dot \tau^m(s)+\gamma(x(s))^m_{kn}\tau^k(s)\tau^n(s)=0,
\end{equation}
\noindent where
\begin{equation}\label{eq32}
\gamma^m_{kn}:=\partial^2_{kn'}\varphi^m(x,t,t')|_{t=t'=0}=
\partial_k\lambda(x,t)^m_n|_{t=0}=\partial_n\mu(x,t)^m_k|_{t=0},
\end{equation}
\noindent Insertion of expression (3) into the first equality of
formula (32), with consideration of definition (8) results in
$\gamma^m_{kn}=h^m_\mu(\Gamma^\mu_{k\nu}h^k_kh^\nu_n+h^k_k\partial_kh^\mu_n)$,
accordingly, functions $\gamma^m_{kn}$ are affine connection
coefficients, while equation (31) is the equation of geodesic in
affine vielbein.

Differentiating equation (30) with respect to $s_2$  and
subsequently assuming $s=0$  and substituting $s_2$  for $s$ ,
with the use of equation (31) we derive the equation for auxiliary
matrices $\mu(x,t)^m_n$:
\begin{equation}\label{eq33}
[\partial_k\mu(x,t)^m_n-\gamma(x)^m_{kn}]\tau^{'k}\tau^{'n}=0,
\end{equation}
\noindent analogous to equation (18). Here according to (26) it is
taken that  $t(0,s)=s\tau$. Insertion of expression (6) of
matrices $\mu(x,t)^m_n$ through deformation functions
$H^m(x,\tilde t)$ into equation (33) once again leads to equation
(19) which provides for canonicity of deformations. Thus condition
of deformations canonicity can be applied not only as equation
(14) or equation (29) equivalent to it, but also as equation (30).

Consequent differentiation of equation (33) with respect to  $s$
in zero leads to the following condition for auxiliary matrices
$\mu(x,t)^m_n$ of canonic deformed group  $\Gamma^H_T$:
\[ \partial^{(n)}_{\{k_1,...,k_n}\mu(x,t)^m_{k_{n+1}\}}|_{t=0}=0,
\]
\noindent in particular
\[ \partial^{(n)}_{\{kl}\mu(x,t)^m_{s\}}|_{t=0}=\rho^m_{\{kls\}}=0
\]
\noindent wherefrom directly follows expression (23) of
coefficients $\rho^m_{kls}$  through curvature tensor $R^m_{kls}$.

So, we have proved the following.

\textbf{Proposition 4 (criterion of multiplication law
canonicity).}

\emph{Deformed group of diffeomorphisms $\Gamma^H_T$  possesses
canonic multiplication law (28) along geodesics if and only if it
is canonic, for auxiliary matrices $\lambda(x,t)^m$  and
$\mu(x,t)^m$ of group $\Gamma^H_T$ equations (29) and (30) being
fulfilled respectively.}

\section{Canonicity of Riemannian space deformations}

In work [1] with requirement for vector length conservation during
finite parallel transports for deformation functions in coordinate
basis we derived equation (11) which uniquely defines them with
respect to metric tensor  $g(x)_{\mu\nu}$ including coefficients
$\Gamma(x)^\sigma_{\mu\nu}$ in the second order of expansion with
respect to $\tilde t$ , which appear to be equal to Cristoffel
symbols (12).

The canonicity of such deformations follows from the fact that
equation (11) written down relatively to functions
$\lambda(x,x'-x)^\mu_\sigma$:
\[\lambda(x,x'-x)^\rho_\mu\lambda(x,x'-x)^\sigma_\nu g(x)_{\rho\sigma}=g(x')_{\mu\nu}
\]
\noindent provides for fulfillment of composition law (24) for
them along geodesics, and, therefore, according to Proposition 2,
such functions correspond to canonic groups  $\Gamma^H_T$. This
can be directly verified if equation (11) is differentiated with
respect to  $\tilde t$ and the result
\[ \partial_{\tilde \lambda} g(x+\tilde t)_{\mu\nu}=2\partial^2_{\tilde \lambda\{ \tilde
\mu}H^\rho(x, \tilde t)\partial_{\tilde \nu\}}H^\sigma(x, \tilde
t)g(x)_{\rho\sigma}
\]
\noindent is inserted into expression (12), taken in point
$x+\tilde t$ . After cancellations we derive:
\begin{equation}\label{eq34}
\Gamma(x+\tilde t)^\sigma_{\mu\nu}g(x+\tilde
t)_{\sigma\lambda}=\partial^2_{\tilde \mu\ \tilde \nu}H^\rho(x,
\tilde t)\partial_{\tilde \lambda}H^\sigma(x, \tilde
t)g(x)_{\rho\sigma}
\end{equation}
\noindent Using once again equation (11) for the expression of
tensor  $g(x+\tilde t)_{\sigma\lambda}$ through derivatives from
deformation functions and considering invertability of matrices
$\partial_{\tilde \lambda}H^\sigma(x, \tilde t)$ we find that
equation (19) follows from equation (34), which provides for
canonicity of transformations, functions of which satisfy the
equation (11).

Let's point out that in Riemannian space geodesics in natural
parameterization are extremals of energy functional
\begin{equation}\label{eq35}
S(x,x')=\frac{1}{2}\int_0^s
g(x(\alpha))_{\mu\nu}\tau^\mu(\alpha)\tau^\mu(\alpha)d\alpha.
\end{equation}
\noindent The function of central field of extremals which start
from each point $x$, if the function of action $S(x,x')$ which is
defined by integral (35) with condition its extremality, satisfied
the Hamilton-Jacobi equation [4]:
\begin{equation}\label{eq36}
g(x')^{\mu\nu}\partial_{\mu'}S(x,x')\partial_{\nu'}S(x,x')=g(x)_{\mu\nu}\tau^\mu\tau^\nu
\end{equation}
\noindent and allows to find vectors  $\tau'$ tangent to geodesic
in current point $x'$  according to the formula:
\[ \tau^{'\nu}=g(x')^{\mu\nu}\partial_{\mu'}S(x,x')
\]
\noindent On the other hand, this vector may be found by means of
parallel transport of initial vector $\tau$  to the point  $x'$
which yields:
\[\partial_{\mu'}S(x,x')=g(x')_{\mu\sigma}\partial_\nu
H^\sigma(x',x-x')\tau^\nu.
\]
\noindent Insertion of this relation into Hamilton-Jacobi equation
(36) leads to fulfillment of equation (11) along geodesic, which
can also be called Hamilton-Jacobi equation for deformation
functions.

Thus, we have proved the following.

\vskip 50mm

\textbf{Proposition 5 (of canonicity of Riemannian space
deformations).}

\emph{Deformed groups of diffeomorphisms  $\Gamma^H_T$ which
specify by their action in $O$ the structure of Riemannian space,
derived with the help of deformations, functions of which satisfy
Hamilton-Jacobi equation (11), are canonic.}

To crown it all, let's point out that finite parallel transports
are naturally joined in the so called groups of parallel
transports $DT=\Gamma^H_T\times)GL^g(n)$
($DT=\Gamma^H_T\times)SO^g(n)$ in case of Riemannian space) [6],
which acts in tangent bundle of space $O$ and is a fundamental
group (in terms of Klein's Erlangen Program) of space with affine
connection (respectively to Riemannian space), while in case of
canonic group $\Gamma^H_T$ group $DT$ has in the set of pure
translations without additional rotations (which in case of curved
space does not form subgroups of group $DT$) one-parameter
subgroups, which translate two arbitrary points $x,x'\in O$  one
into another. On the contrary, the existence of such subgroups in
group of parallel transports $DT$ fulfills canonicity of group
$\Gamma^H_T$ .


\end{document}